\documentclass[12pt]{article}
\usepackage{amssymb,amsbsy,amsmath,amsfonts,amssymb,amscd, graphics}
\usepackage{latexsym,exscale,epsfig}
\usepackage[mathscr]{euscript}

\newcommand{\noi}{\noindent}

\newcommand{\RR}{\mathbb R}

\newcommand{\ZZ}{{\mathbb Z}}
\renewcommand{\epsilon}{\varepsilon}
\newcommand{\pf}{{\mathfrak p}_F}
\newcommand{\of}{{\mathfrak o}_F}
\newcommand{\ofr}{\mathfrak o}
\newcommand{\pfr}{\mathfrak p}
\newcommand{\Gfr}{\mathfrak g}
\newcommand{\tGfr}{\tilde \Gfr}
\newcommand{\Hfr}{\mathfrak h}
\newcommand{\tHfr}{\tilde \Hfr}
\begin{document}

\setcounter{section}{0}

\vskip1cm

\centerline{\bf Buildings of classical groups and centralizers}

\centerline{\bf of Lie algebra elements}

\bigskip

\centerline{\sc P. Broussous and S. Stevens}
\medskip

\centerline{\it F\'evrier 2004}

\vskip1cm

\centerline{\bf Abstract}
\medskip

Let $F_o$ be a non-archimedean locally compact field of residual 
characteristic not $2$. Let $G$ be a classical group over $F_o$
(with no quaternionic algebra involved) which is not of type $A_n$ for
$n>1$. Let $\beta$ be an element of the Lie algebra $\Gfr$ of $G$
that we assume semisimple for
simplicity. Let $H$ be the centralizer of $\beta$ in $G$ and $\Hfr$
its Lie algebra. Let $I$ and $I_{\beta}^{1}$ denote the (enlarged)
Bruhat-Tits buildings of $G$ and $H$ respectively. We prove that
there is a natural set of maps $j_{\beta}$ :
$I_{\beta}^{1}\rightarrow I$ which enjoy the following properties:
they are affine, $H$-equivariant, map any apartment of
$I_{\beta}^{1}$ into an apartment of $I$ and are compatible with the
Lie algebra filtrations of $\Gfr$ and $\Hfr$. In a particular case,
where this set is reduced to one element, we prove that $j_{\beta}$
is characterized by the last property in the list. We also prove a
similar characterization result for the general linear group.
\vskip1cm

%
\noi {\bf Introduction}
\smallskip

In this paper we establish new functoriality properties between affine
Bruhat-Tits  buildings of classical reductive groups over local
fields. More precisely let $F_o$ be a non-archimedean local field of
residual characteristic not $2$ and
$G$ be the group of $F_o$-rational points of a classical group defined
over $F_o$. We assume that $G$ is the isometry group of an
$\epsilon$-hermitian form over an $F$-vector space, where $F$ is a
(commutative) extension of $F_o$ of degree less than $2$.  We denote
by $\mathfrak g$ the Lie algebra of $G$ and by
$I$ its affine building. Let $\beta$ be an element of $\mathfrak g$
that we assume to be semisimple for simplicity.  Let $H$ be the
centralizer of $\beta$ in $G$. Then $H$ is the group of $F_o$-rational
points of a product of groups of the
form ${\rm Res}_{E_o/F_o}{\boldsymbol H}_i$, where $E_o /F_o$ is a
field extension 
and where ${\rm Res}$ denotes
the functor of restriction of scalars. Here $i$ runs over a finite set
$\tilde J$. Each ${\boldsymbol H}_i$ is either a classical group as above or a general
linear group. We denote by $J_+\subset \tilde J$ the (possibly empty) subset
of indices corresponding to linear groups. We denote by $\mathfrak h$
the Lie algebra of $H$ and by $I_{\beta}^{1}$ its (enlarged) affine
building. Then there is a natural set of maps $j_{\beta}~:$ 
$I_{\beta}^{1}\rightarrow I$ which depend on  identifications of the
enlarged buildings of ${\boldsymbol H}_i$, $i\in J_+$, with certain sets of lattice
functions (see {\S}4 below). In particular, when $J_+ =\emptyset$, there is a natural
choice of $j_\beta$. The maps $j_{\beta}$ enjoys the following
properties:
\smallskip

\noi a) They are  affine.

\noi b) They are $H$-equivariant.

\noi c) They map any apartment of $I_{\beta}^{1}$ into an apartment of $I$.

\noi d) They are compatible with the Lie algebra filtrations of
$\mathfrak g$ and $\mathfrak h$ (cf. {\S}9). 
\smallskip

In \cite{BL} it was proved that when $G$ is the general linear
group and $\beta$ is an elliptic element then, replacing the buildings
by the non-enlarged buildings, there is such a natural
map $j_{\beta}$ satisfying the conditions above. It is actually
characterized  by properties a) and b). However in the case of a
classical group (and assuming that $J_+ =\emptyset$) it is no longer
 true that properties a), b) and c) characterize $j_{\beta}$. The
simplest counter-example is the following. Consider the case of
$G={\rm Sp}_{2}(F_o)={\rm SL}(2,F_o)$. One may choose $\beta$ is such
a way that $H$ is $E^1$, the group of norm $1$ elements of a ramified
quadratic extension $E /F_o$. Then $I_{\beta}^{1}$ is reduced to a point
and fixing a map $j_{\beta}$ satisfying a) b) and c) amounts to choosing a point in
$I$ fixed by the torus $E^{1}$. But $E^{1}$ is contained in an
Iwahori subgroup of $G$ and therefore fixes a chamber of $I$.

We prove that, in the case of a general linear group and of an elliptic element
$\beta$, the map $j_{\beta}$ of \cite{BL} is actually characterized
by property d). In the case of a classical group, we also prove that
if $J_+ =\emptyset$ and if a technical condition on $\beta$ is
satisfied  then $j_{\beta}$ is characterized by condition d). We
conjecture that when $J_+ =\emptyset$ then $j_{\beta}$ is indeed
characterized by property d).  

\medskip

In this work, we do not actually assume $\beta$ to be semisimple but only to satisfy
a weaker assumption (see hypothesis (H1) of {\S}5). Such elements
naturally appear in the generalization of the theory of strata due to
Bushnell and Kutzko \cite{BK} to the case of classical groups (see
the work of the second author \cite{S1}, \cite{S2}. Even though the work of the second author does not  use
the  theory of affine buildings  in a straightforward way (it uses the equivalent
language of hereditary orders), the existence and properties of the maps $j_{\beta}$
are applied to the representation theory of $G$, particularly in \cite{S2}.    
\medskip

The paper is organized as follows. In {\S}2 we recall the structure
of the maximal split tori of $G$. In {\S}3,4, using ideas of Bruhat and
Tits, we give a model of the affine building of $G$ in terms of
``self-dual lattice functions''. In {\S}5 we study the centralizers
in $\mathfrak g$ and $G$ of the Lie algebra element $\beta$. The
construction of the maps $j_{\beta}$ is done in {\S}6 and their
properties are established in {\S}7,8 and 9. In {\S}10 We prove  the
uniqueness result for the general linear group and finally {\S}11 is
devoted to the uniqueness result in the classical group case. 
\medskip

 We thank G. Henniart, B. Lemaire, G. Prasad and J.-K. Yu for
 stimulating discussions.

\noi {\bf 1. Notation}
\smallskip 

Here $F_o$ is the ground field; it is assumed to be non-archimedean,
locally compact and equipped with a discrete valuation $v$ normalized
in such a way that $v(F_{o}^{\times})$ is the additive group of
integers.  We assume that the residual characteristic of $F_o$ is not
$2$. We fix a Galois extension $F/F_o$ such that $[F:F_o]\leqslant 2$
and set $\sigma_F = {\rm id}_{F}$ if $F=F_o$ and take $\sigma_F$ to be
the generator of ${\rm Gal}(F/F_o )$ in the other case. We still
denote by $v$ the unique extension of $v$ to $F$. We fix 
$\epsilon\in \{ \pm 1\}$ and a finite dimensional left $F$-vector 
space $V$. Recall that a $\sigma_F$-skew form $h$ on $V$ is a 
$\ZZ$-bilinear map $V\times V\rightarrow F$ such that
$$
h(\lambda x ,\mu y)= \lambda^{\sigma_F}\mu h(x,y)\ , \ \lambda ,\mu\in
F,\  x , y\in V \ .
$$

\noi Such a form is called $\epsilon$-hermitian if $h(y,x)=\epsilon
h(x,y)^{\sigma_F}$ for all $x$, $y\in V$. From now on we fix such an
$\epsilon$-hermitian form on $V$ and we assume it is non-degenerate
(the orthogonal of $V$ is $\{ 0\}$).

For $a\in {\rm End}_{F}(V)$, we denote by $a^{\sigma_h}= a^{\sigma}$
 the  adjoint of $a$ with
respect to $h$, i.e. the unique $F$-endomorphism of $V$ satisfying 
$h(ax,y) = h(x,a^{\sigma}y)$ for all $x$, $y\in V$.

\smallskip

 We denote by $\boldsymbol G$ the simple algebraic $F_o$-group
 whose set of $F_o$-rational points $G$ is formed of the $g\in {\rm
 GL}_{F}(V)$ satisfying $g.h = h$ (it is not necessarily connected). Here $g.h$ is
 the  form given by
 $g.h(x,y)= h(gx,gy)$, $x$, $y\in V$.

\smallskip

 We know (\cite{Sch}(6.6), page 260) that in the case $\sigma_F \not=
 {\rm id}_F$, we may reduce to the case $\epsilon =1$. So we have three
 possibilities:
\smallskip

$\sigma_F ={\rm id}_F$ and $\epsilon =1$, the orthogonal case;

$\sigma_F ={\rm id}_F$ and $\epsilon =-1$, the symplectic case;

$\sigma_F \not= {\rm id}_F$ and $\epsilon =1$, the unitary case.

\smallskip

 We abbreviate ${\tilde G} = {\rm GL}_F (V)$ and  $\tGfr = {\rm
 End}_F (V)$.

\medskip

\noi {\bf 2. The maximal split tori of $\boldsymbol G$}
\smallskip

Recall that a subspace  $W\subset V$ is totally isotropic if
$h(W,W)=0$ and that maximal such subspaces have the same dimension
$r$, the Witt index of $h$. Set $I=\{\pm 1 , \pm 2 , \dots,\pm
r\}$ and $I_o =\{ (0,k)\ ; \ k=1,\dots ,n-2r\}$. We fix a {\it Witt
decomposition} of $V$, that is
\smallskip

\noi  two maximal totally isotropic subspaces
$V_+$ and $V_-$,

\noi bases $(e_i )_{i=1,\dots ,r}$, $(e_{-i})_{i=1,\dots ,r}$, $(e_i
)_{i\in I_o}$ of $V_+$, $V_-$ and $V_o :=(V_+ +V_- )^{\perp}$,
\smallskip

\noi such that
\smallskip

\noi $h(e_i ,e_i )=0$, $i\in I$,

\noi $h(e_i ,e_j ) =0$, for $i$, $j\in I$ with $j\not=-i$ or $i\in I$,
$j\in I_o$,

\noi $h(e_i ,e_{-i})=1$, for $i\in I$ with $i>0$,

\noi $h(x,x)\not=0$, for $x\in V_o$ and $x\not= 0$.
\smallskip

The Witt decomposition gives rise to a maximal $F_o$-split torus
$\boldsymbol S$ whose group of $F_o$-rational points is
$$
S=\{ s\in G\ ; \ se_i \in F_o e_i\ , \ i\in I\ {\rm and}\ (s-{\rm
Id})V_o =0\} \ .
$$

\noi It has dimension $r$, the $F_o$-rank of $\boldsymbol
G$. Conversely any maximal $F_o$-split torus of $\boldsymbol G$ is
obtained from a Witt decomposition as above.  The
centralizer $\boldsymbol Z$ of $\boldsymbol S$ in $\boldsymbol G$ has for
$F_o$-rational points
$$
Z=\{ z\in G\ ; \ ze_i \in Fe_i\ , \ i\in I\ {\rm and }\ zV_o =V_o\} \
.
$$

For each $i\in I$, we have a morphism of algebraic groups $a_i$~:
$\boldsymbol Z \rightarrow {\rm Res}_{F/F_o}({\boldsymbol G}_m)$ given by $ze_i =a_i
(z)e_i$. Note that $a_{-i}(z) = a_{i}(z)^{-\sigma}$. We also denote by
$a_i$~: ${\boldsymbol S}\rightarrow {\boldsymbol G}_{m\  /F_o}$ the character
obtained by restriction. We have $a_i =-a_{-i}$ in $X^{*}({\boldsymbol
S})$, the $\mathbb Z$-module of rational characters of $\boldsymbol
S$. The $a_i$, $i\in I$, $i>0$, form a basis of  $X^{*}({\boldsymbol
S})$. 
\smallskip

The normalizer $\boldsymbol N$ of $\boldsymbol Z$ in $\boldsymbol G$ is the
sub-algebraic group whose $F_o$-rational points are the elements of
$G$ which stabilize $X_o$ and permute the lines $V_i = Fe_i$, $i\in
I$. The group $N={\boldsymbol N}(F_o )$ is the semidirect product of $Z$
by the subgroup $N'$ formed of the elements which permute the $\pm
e_i$, $i\in I$.
\smallskip

\medskip

\noi {\bf 3.  $\rm MM$-norms and self-dual lattice-functions}

\medskip
 We keep the notation as in the previous sections.
\smallskip

 Recall that a {\it norm} on $V$ is a map $\alpha$~: $V\rightarrow
 {\mathbb R}\cup \{\infty\}$ satisfying:
\smallskip

i) $\alpha (x+y)\geqslant {\rm Inf}(\alpha (x),\alpha (y))$, $x$, $y\in V$,

ii) $\alpha (\lambda x)= v(\lambda )+\alpha (x)$, $\lambda\in F$,
$x\in V$,

iii) $\alpha (x)=\infty$ if and only if $x=0$.
\smallskip

 We denote by ${\rm Norm}^1 (V)$ the set of norms on $V$.
\smallskip

\noi {\bf  Definition 3.1.} (cf. [BT](2.1)) {\it Let $\alpha\in
{\rm Norm}^1 (V)$. We say that $\alpha$ is dominated by $h$ if 
$$
\alpha (x)+\alpha (y)\leqslant v(h(x,y))\ {\rm for \ all}\ x,\, y\in V\
.
$$

\noi We say that $\alpha$ is an $\rm MM$-norm for $h$ 
(\/}maximinorante {\it in french), if $\alpha$ is a maximal element of the
set of norms dominated by $h$.}

In [BT](2.5) an involution $\bar{}$ is defined on ${\rm Norm}^1 (V)$
in the following way. If $\alpha\in {\rm Norm}^1 (V)$, then 
$$
\bar{\alpha}(x)=\inf_{y\in V}[v(h(x,y))-\alpha (y)]\ , \ x\in V\ .
$$

\noi We then have
\smallskip

\noi {\bf  Proposition 3.2.} (cf. [BT](Prop.\ 2.5)) {\it An element
$\alpha$ of ${\rm Norm}^1 (V)$ is an $\rm MM$-norm if and only
 if $\bar{\alpha}=\alpha$.}
\smallskip

We are going to describe the set ${\rm Norm}^{1}_{h}(V)$  of   $MM$-norms
 in terms of self-dual
lattice-functions. Recall \cite{BL} that a lattice-function in $V$ is a
 function $\Lambda$ which maps a real number to an $\of$-lattice in
$V$ and satisfies:
\smallskip

i) $\Lambda (r)\subset \Lambda (s)$ for $r\geqslant s$, $r$, $s\in
\mathbb R$,

ii) $\Lambda (r+v(\pi_F )) = \pf \Lambda (r)$, $r\in \mathbb R$,

iii) $\Lambda$ is left-continuous.
\smallskip

\noi Here $\of$ denotes the ring of integers of $F$, $\pf$ the maximal
ideal of $\of$ and $\pi_F$ a uniformizer of $F$. As in \cite{BL}, we denote
by ${\rm Latt}_{\ofr_F}^{1} (V)$ (or by ${\rm Latt}^{1}(V)$ when no
confusion may occur)  the set of $\of$-lattice-functions in $V$.
\smallskip
 Recall \cite{BL} that ${\rm Norm}^1 (V)$ and ${\rm Latt}^1 (V)$ may be
 canonically identified in the following way. To $\alpha\in{\rm
 Norm}^1 (V)$, we attach the function $\Lambda = \Lambda_{\alpha}$
 given by
$$
\Lambda (r) =\{ x\in V \ ; \ \alpha (x)\geqslant r\}\ , \ r\in
{\mathbb R}\ .
$$

\noi Conversely a lattice-function $\Lambda$ corresponds to the norm
$\alpha$ given by
$$
\alpha (x) = \sup\{ r\ ; \ x\in \Lambda (r)\}\ , \ x\in V\ .
$$

\noi  For a $\Lambda\in {\rm Latt}^1 (V)$ and $r\in \RR$, set 
$$
\Lambda (r+)=\bigcup_{s>r}\Lambda (s)\ .
$$
\noi For an $\of$-lattice $L$ in $V$, we define its dual $L^{\sharp} = 
L^{\sharp_h}$  by 
$$
L^{\sharp} = \{ x\in V\ ; \ h(x,L)\subset \pf\} \ .
$$

\noi Finally, we define the dual $\Lambda^{\sharp} = \Lambda^{\sharp_h}$ of a
 lattice-function  $\Lambda$ by 
$$
\Lambda^{\sharp}(r) = [\Lambda ((-r)+)]^{\sharp} \ , r\in \RR\ .
$$  

\noi We say that a lattice function $\Lambda$ is self dual if 
$\Lambda^{\sharp}=\Lambda$ and we denote by $Latt^{1}_{h}(V)$ the
 corresponding  set. 
\smallskip

\noi {\bf  Proposition 3.3}. {\it Given a norm $\alpha\in {\rm Norm}^1
(V)$, we have $\Lambda_{\bar\alpha} = \Lambda_{\alpha}^{\sharp}$.}
\smallskip

\noi {\bf  Corollary 3.4}. {\it Let $\alpha$ be a norm on $V$. Then
$\alpha$ is an MM-norm if and only if the attached lattice-function
$\Lambda$ is self-dual.}
\smallskip

\noi {\it Proof of Proposition}. Let $x\in V$ and $r\in \RR$. Then the
fact that $x\in \Lambda_{\bar \alpha}(r)\backslash \Lambda_{\bar
\alpha}(r+)$ is equivalent to the following points:
\smallskip

\noi i) ${\bar \alpha}(x)=r$;

\noi ii) there exists $y\in V$ such that $v(h(x,y))-\alpha (y) = r$, and
for all $y\in V$, we have $v(h(x,y))-\alpha (y)\geqslant r$;

\noi iii) there exists  $y\in V$ such that $v(h(x,y))=0$ and $\alpha
(y)=-r$, and for all $y\in V$ such that $\alpha (y)>-r$, we have
$v(h(x,y))>0$ (scale by a suitable power of  a uniformizer $\pi_F$);

\noi iv) there exists $y\in \Lambda_{\alpha}(-r)\backslash
\Lambda_{\alpha}(-r+)$ such that $h(x,y)\in \of\backslash\pf$, and for
all $y\in \Lambda_{\alpha}(-r+)$ we have $h(x,y)\in \pf$;

\noi v) $x\in
\Lambda_{\alpha}^{\sharp}(r)\backslash\Lambda_{\alpha}^{\sharp}(r+)$. 
\smallskip

\noi This proves  that the two lattice-functions $\Lambda_{\bar
\alpha}$ and $\Lambda_{\alpha}^{\sharp}$ share the same discontinuity
points and that at those points they take the same values; so there are
equal. 
 
\medskip

 Let ${\rm Norm}^2 \tGfr $ (resp. ${\rm Latt}^2 \Gfr$) denote the
 $\tilde G$-set of square norms in $\tGfr$ (resp. of square
 lattice-functions in $\tGfr$; see \cite{BT1} and \cite{BL}). Recall that a
 lattice-function $\Lambda^2$ in the $F$-vector space $\tGfr$ is
 square if there exists $\Lambda\in {\rm Latt}^1 (V)$ such that
 $\Lambda^2 = {\rm End}(\Lambda )$, where 
$$
{\rm End}(\Lambda )(r) = \{a\in \tGfr\ ; \ a\Lambda (s)\subset \Lambda
(s+r),\ s\in {\mathbb R}\}, \ r\in {\mathbb R}\ .
$$

\noi An additive  norm on $\tGfr$ is square if the corresponding
lattice function is square. Recall \cite{BT1} that  ${\rm Norm}^{1}(V)$
and ${\rm Norm}^2 \tGfr$
(and therefore ${\rm Latt}^{1}(V)$ and ${\rm Latt}^2 \tGfr$ by
transfer of structure) are  endowed with  affine
structures~: the barycenter of two points with positive  weights is
defined. 

 The involution $\sigma$ acts on ${\rm Norm}^2 \tGfr$ via
$$
\alpha^{\sigma}(a)=\alpha (a^{\sigma}),\ a\in \tGfr , \ \alpha\in {\rm
Norm}^2 \tGfr\ .
$$
\noi By transfer of structure, $\sigma$ acts on ${\rm Latt}^2 \tGfr$
via
$$
\Lambda^{\sigma}(r)=[\Lambda (r)]^{\sigma} , \ \Lambda\in {\rm Latt}^2
\tGfr , r\in {\mathbb R}\ .
$$

\noi A square norm $\alpha$ (resp. a square lattice function
$\Lambda$) is said to be self-dual if $\alpha =\alpha^{\sigma}$
(resp. $\Lambda =\Lambda^{\sigma}$). We denote by
 ${\rm Norm}^{2}_{\sigma}\tGfr$ and ${\rm Latt}^{2}_{\sigma}\tGfr$ 
the corresponding sets.  

 Now, in terms of lattice functions,  Corollary 2 of \cite{BT2}, page 163,
 writes:

\smallskip

\noi {\bf Lemma 3.5} {\it The map $\Lambda\mapsto {\rm End}(\Lambda )$
induces a bijection from the set of self-dual lattice functions in $V$
to the set of self-dual square lattice functions in  $\tGfr$.}
\medskip

In other words, for any $\Lambda\in {\rm Latt}^{2}_{\sigma}\tGfr$, there exists
a unique $\Lambda^2 = \Lambda_{h}^{2}\in {\rm Latt}^{1}_{h}(V)$ such that 
${\rm End}(\Lambda )=\Lambda^2$.
\smallskip

  Note that the sets ${\rm Latt}^{1}_{h}(V)$, 
${\rm Norm}^{1}_{h}(V)$, ${\rm Latt}^{2}_{\sigma}\tGfr$ and 
 ${\rm Norm}^{2}_{\sigma}\tGfr$ are $G$-sets and that the
 various identifications among them are $G$-equivariant. 
\smallskip

 Let $u\in F^{\times}$ and assume that $uh$ is still an 
$\epsilon$-hermitian form with respect to $\sigma_F$. Then the involution
$\sigma$ of $\tGfr$ corresponding to $uh$ remains the same and defines
the same unitary group $G\subset {\tilde G}$. For $\Lambda\in {\rm Latt}^{1}(V)$ 
and $s\in \mathbb R$, we denote by $\Lambda +s$ the lattice function given by 
$(\Lambda +s)(r)= \Lambda (s+r)$, $r\in \mathbb R$.  
\smallskip

\noi {\bf Lemma 3.6}. {\it Let $\Lambda^2 \in {\rm Latt}^{2}_{\sigma}\tGfr$ 
and $\Lambda^{2}_{h}$ (resp. $\Lambda_{uh}^{2}$) be the unique element of 
${\rm Latt}^{1}_{h}(V)$ (resp. of ${\rm Latt}^{1}_{uh}(V)$) satisfying 
${\rm End}(\Lambda_{h}^{2})=\Lambda^{2}$ (resp. ${\rm End}(\Lambda_{uh}^{2})
=\Lambda^{2}$). Then $\Lambda_{uh}^{2}= \Lambda_{h}^{2}-v(u)/2$, that is 
$\Lambda^{2}_{uh} (r) = \Lambda^{2}_{h}(r-v(u)/2)$, $r\in \mathbb R$.}
\smallskip

\noi {\it Proof}. We easily check that for $\Lambda\in {\rm Latt}^{1}(V)$ and
$s\in \mathbb R$, we have
$$
\Lambda^{\sharp_{uh}} =u^{-\sigma}\Lambda^{\sharp_{h}}\ {\rm and} \ 
(\Lambda +s)^{\sharp_{h}} = \Lambda -s\ .
$$
We certainly have ${\rm End}(\Lambda_{h}^{2}-v(u)/2)={\rm End}(\Lambda^{2}_{h})$
$=$ $\Lambda^2$. So by a unicity argument, we must prove that $\Lambda_{h}^{2}-v(u)/2\in 
{\rm Latt}_{uh}^{1}(V)$. But 
$$
(\Lambda_{h}^{2}-v(u)/2)^{\sharp_{uh}}=u^{-\sigma}(\Lambda_{h}^{2}-v(u)/2)^{\sharp_{h}}
$$
$$
=u^{-\sigma}(\Lambda_{h}^{2}+v(u)/2) = \Lambda^{2}_{h}+v(u)/2-v(u^{\sigma}) = \Lambda_{h}^{2}-v(u)/2\ ,
$$
\noi as required.

\medskip 

\noi {\bf 4. The building as a set of self-dual lattice-functions}

\medskip

 Let $I$ denote the building of the standard valuated root datum of $G$
 introduced in \cite{BT2} and $A$ denote the apartment of $I$ attached to
 $\boldsymbol S$. Write $V^{*} =X^* ({\boldsymbol S}\otimes {\mathbb R})$; this
 is an $\mathbb R$-vector space with basis $(a_i )_{i=1,\dots
 ,r}$. Let $V$ denote the linear dual of $V^*$. We identify $A$ with $V$.

 To a point $p\in A\simeq V$, we attach the norm $\alpha_p$ on $V$
 defined by
$$
\alpha_p (\sum_{i\in I}\lambda_i e_i +x_o )=\inf [\omega (x_o
),\inf_{i\in I}(v(\lambda_i )-a_i (p))], \ x_o \in V_o ,\
\lambda_i \in F \ {\rm for }\ i\in I\ .
$$
\noi Here $\omega (x_o ) = \frac{1}{2}v(h(x_o ,x_o ))$, $x_o \in V_o$. 
\smallskip

Here are two important facts from \cite{BT2}.

\smallskip

\noi {\bf Proposition 4.1.} (\cite{BT2}(Prop.\ 2.9, 2.11(i))) {\it The map
$p\mapsto \alpha_p$ is a bijection from $A$ to the set of $MM$-norms
on $V$ which split in the decomposition $\displaystyle V=\oplus_{i\in
I}Fe_i \oplus V_o$. It is $N$-equivariant.}
\smallskip

 For the notion of splitting for norms, see \cite{BT1}(1.4). 
\smallskip

\noi {\bf Proposition 4.2.} (\cite{BT2}(2.12)) {\it i) The map $p\mapsto
\alpha_p $ extends in a unique way to a $G$-equivariant and affine
bijection $j_h$~: $I\rightarrow {\rm Norm}^{1}_{h}(V)$ (in particular
${\rm Norm}^{1}_{h}(V)$ is a convex subset of ${\rm Norm}^1 (V)$).

\noi ii) The map $j_h$ is the unique affine and $G$-equivariant map 
$I\rightarrow {\rm Norm}^{1}_{h}(V)$.}
\smallskip

From {\S}3, we get a unique affine and $G$-equivariant map $I\rightarrow
 {\rm Latt}_{h}^{1}(V)$ that we still denote by $j_{h}$.
\smallskip

  For $r\in \mathbb R$, let
${\mathscr V}_{o}^r$ be the lattice of $V_o$ given by $\{ x_o \in V_o
\ ;  \ \omega (x_o )\geqslant r\}$. For $x\in \mathbb R$, let
$\lceil x \rceil$ denote the least integer greater than or equal
to $x$. Then the map $j_{h}~:$ $I\rightarrow {\rm Latt}_{h}^{1}(V)$ 
  is given on $A$ by $j_{h}(p) = \Lambda_p$, where
$$
\Lambda_p (r) = {\mathscr V }_{o}^{r}\oplus\bigoplus_{i\in
I}{\mathfrak p}_{F}^{\lceil r+a_i (p)\rceil}e_i \ , \ r\in {\mathbb
R}\ .
$$

\smallskip

Let $u$ be an element of $F^{\times}$ such that $uh$ remains 
$\epsilon$-hermitian with respect to $\sigma_F$. It follows from the proof 
of Lemma (3.6) that if $\Lambda\in {\rm Latt}^{1}(V)$, we have $\Lambda\in 
{\rm Latt}^{1}_{h}(V)$ if, and only if, $\Lambda -v(u)/2\in 
{\rm Latt}^{1}_{uh}(V)$. Since ${\rm End}(\Lambda +s)={\rm
End}(\Lambda )$, for  $\Lambda\in {\rm Latt}^{1}(V)$ and $s\in \RR$, 
 the bijective map $j_{\sigma}~:$ 
$I \rightarrow {\rm Latt}^{2}_{\sigma}(V)$, given by 
$j_{\sigma}={\rm End}\circ j_{h}$, does not depend on the choice of the 
form $h$, the involution $\sigma$ being fixed. By construction it is affine 
and $G$-equivariant. It is uniquely determined by these two properties. 
Indeed if $j_{\sigma}'~:$ $I\rightarrow {\rm Latt}^{2}_{\sigma}(V)$ is affine 
and $G$-equivariant, so is $(j_{\sigma}')^{-1}\circ j_{\sigma}~:$ 
$I\rightarrow I$. But such a map must be the identity map. 

\medskip

We also recall here the description of the enlarged building $I^1$ 
of $\tilde G = {\rm GL}_F(V)$ in terms of lattice functions.  

\smallskip

\noi {\bf Proposition 4.3.} (\cite{BT1}(2.11)) {\it i) There is a
 $\tilde G$-equivariant and affine
bijection $j$~: $I^1\rightarrow {\rm Norm}^{1}(V)$.

\noi ii) If we have another affine and $\tilde G$-equivariant map 
$j':I^1\rightarrow {\rm Norm}^{1}(V)$ then there exists 
$r\in\mathbb R$ such that, for all $\alpha\in {\rm Norm}^{1}(V)$, 
$j'(\alpha)=j(\alpha)+r$.}
\smallskip

From \cite{BL} Proposition 2.4, for each $j$ as in Proposition 4.3, we get an 
affine and $\tilde G$-equivariant map $I^1\rightarrow {\rm Latt}^{1}(V)$ 
that we also denote by $j$.

\bigskip

\noi {\bf 5. Centralizers of Lie algebra elements}
\medskip

We denote by $\Gfr$ the Lie algebra of  $G$:
$$
\Gfr =\{ a\in \tGfr\ ; \ a+a^{\sigma}=0\} \ .
$$
We consider an element $\beta$ of $\Gfr$ satisfying
$$
{\rm The \ }F\hbox{\rm-algebra}\ E:=F[\beta ]\subset \tGfr\  {\rm is\ a \ direct \
sum \ of \ fields.}
\leqno (\hbox{H1})
$$
We write $\tHfr$ (resp. $\Hfr$) for the centralizer of $\beta$ in $\tGfr$ 
(resp. in $\Gfr$) and ${\tilde H}$ (resp. $H$) for the fixator of $\beta$ in 
$\tilde G$ (resp. in $G$) for the adjoint action. 

\smallskip

Since $\sigma (\beta )=-\beta$, we have easily that $E\subset \tGfr$
is $\sigma$-stable. We write 
$$
E=\bigoplus_{i=1,\dots,t} (E_i\oplus E_{-i})\ \oplus\ 
\bigoplus_{k=1,\dots ,s}E_{(0,k)},
$$
where, for each $i$ in $J=\{\pm 1,\ldots,\pm t\}$ or 
$J_o=\{(0,k):k=1,\ldots,s\}$, $E_i$ is a field extension of $F$, and we have 
labeled the components such that, for each $i\in J_o\cup J$,
$$
\sigma(E_i)=E_{-i},  \leqno (\hbox{H2})
$$
with the understanding that $i=-i$, for $i\in J_o$. We remark that the torus 
$E\cap G$ in $G$ is anisotropic (modulo the centre) if and only if 
$J=\emptyset$ and that every maximal anisotropic torus in $G$ takes this form 
(see \cite{Mor} Proposition 1.3).

For each $i\in J_o$, we set $E_{i}^o = \{ a\in E_i\ ; a=a^{\sigma}\}$, so
that $E_i/E_{i}^o$ is a Galois extension of degree $\leqslant 2$ and
a generator  of ${\rm Gal }(E_i /E_{i}^o )$ is $\sigma_{E_i} :=
\sigma_{\vert E_i }$. For $i\in J_o\cup J$, let ${\mathbf 1}_i$ be the
idempotent of $E$ attached to $E_i$; from (H2), we have $\sigma
({\mathbf 1}_i )= {\mathbf 1}_{-i}$. We have the decomposition
$$
V= \bigoplus_{i\in J_o\cup J} V_i\ , V_i ={\mathbf 1}_i V\ .
$$
\noi Note that, if $i\not= -k$, $v\in V_i$ and $w\in V_k$, we have 
$h(v,w)$ $=$ $h({\mathbf 1}_i v,w)$ $=$ $h(v, {\mathbf 1}_i w)=0$ so, for 
$i\in J_o\cup J$,
$$
V_i^\perp=\bigoplus_{k\ne -i} V_k.
$$

For $i\in J_o\cup J$, $V_i$ is naturally an $E_i$-vector space and
we have obvious isomorphisms of algebras and groups respectively:
$$
\tHfr \simeq \prod_{i\in J_o\cup J} {\rm End}_{E_i}V_i\ ,
$$
$$
{\tilde H} \simeq \prod_{i\in J_o\cup J} {\rm Aut}_{E_i}V_i \ .
$$
\noi The involution $\sigma$ stabilizes $\tHfr \subset \tGfr$ and, for each 
$i$, $\sigma({\rm End}_{E_i}V_i)={\rm End}_{E_{-i}}V_{-i}$. For $i\in J_o$, we 
write $\sigma_i = \sigma_{\vert {\rm End}_{E_i}V_i}$. Let us fix $i\in J_o$. 
The map $\sigma_i$ is an involution of the  central simple $E_i$-algebra ${\rm
End}_{E_i}V_i$. By a classical theorem (\cite{Inv} Theorem 4.2), there exists 
$\epsilon_i\in\{ \pm 1\}$ and a non-degenerate $\epsilon_i$-hermitian form $h_i$
on $V_i$ relative to $\sigma_{E_i}$ such that $\sigma_i$ is the
involution attached to $h_i$. Of course $h_i$ is only defined up to a
scalar in $E_{i}^{\times}$. Let
$$
H_i = \{ g\in  {\rm Aut}_{E_i}V_i\ ; \ gg^{\sigma_i}=1\}
$$
\noi be the unitary group attached to $h_i$. On the other hand, for $i\in J$, 
we put
$$
H_i = {\rm Aut}_{E_i}V_i,
$$
so that $\sigma(H_i)=H_{-i}$ and $H_i$ is isomorphic to $\{g\in H_i\times 
H_{-i}:  gg^{\sigma}=1\}$ by $h\mapsto (h,h^{-\sigma})$. Then, putting 
$J_+=\{1,\ldots,t\}$, we have a natural group isomorphism
$$
H\simeq \prod_{i\in J_o\cup J_+} H_i\ .
$$

We may actually require a compatibility relation between the forms
$h_i$, $i\in J_o$ and the form $h$. Let us fix $i\in J_o$. 
Let $\lambda_i~: E_i\rightarrow F$
be any  $\sigma$-equivariant non-zero $F$-linear form. Such forms
exist. Indeed choose a non-zero linear form $\lambda_i^o~:
E_i^{o}\rightarrow F_o$. If $F=F_o$ then we put $\lambda =\lambda_i^o \circ
{\rm Tr}_{E/E_i^{o}}$. Otherwise $E_i=FE_i^{o}$ and we can extend
$\lambda_i^o$ by linearity to get the required map $\lambda_i$. In
all cases we have:
$$
\lambda_i^o \circ {\rm Tr}_{E_i /E_i^o} = {\rm Tr}_{F/F_o}\circ
\lambda\ .
\leqno (5.1)
$$

 We still write $h$ for the restriction of $h$ to $V_i$.
\smallskip

\noi {\bf Lemma 5.2.} {\it Let $i\in J_o$. There exists a unique 
$\epsilon$-hermitian form $h_i~: V_i \times V_i \rightarrow E_i$ relative to 
$\sigma_{E_i}$ such that
$$
h(v,w)= \lambda_i (h_i (v,w)),\ \hbox{ for\ all}\ v,w\in V_i\ . 
\leqno (5.3)
$$
\noi It is non-degenerate.}
\smallskip

\noi {\it Proof}. Since we have the orthogonal decomposition 
$$
V=V_i \perp \bigoplus_{k\not= i}V_k\ ,
$$
\noi the restriction $h_{\vert V_i}$ is non-degenerate.

 The $F$-linear map ${\rm Hom}_{E_i}(V_i ,E_i
)\rightarrow {\rm Hom}_F (V_i ,F)$, $\varphi\mapsto \lambda_i\circ
\varphi$ is an isomorphism of $F$-vector space. Indeed if $\varphi$
lies in the kernel, we have ${\rm Im}(\varphi )\subset {\rm
Ker}(\lambda_i )$, a strict subspace of $E_i$, and $\varphi$ must be
trivial.  Moreover the two dual spaces have the same
$F$-dimension. For $v\in V_i$ let $h_v$ be the element of 
${\rm Hom}_F (V_i ,F)$ given by $h_v (w) = h(v,w)$. There exists a
unique $\varphi_w\in {\rm Hom}_{E_i}(V_i ,E_i )$ such that $h_v
=\lambda_i \circ \varphi_w$. It is now routine to check that $h_i
(v,w):= \varphi_{v}(w)$, $v,w\in V_i$, has the required properties.

\smallskip

 We easily check that if $h_i$ is as in the lemma, then  the
 corresponding involution on ${\rm End}_{E_i}V_i$ is $\sigma_i$. In
 the following we assume that the forms $h_i$, $i\in J_o$, satisfy
 (5.3).

\smallskip

 For technical reasons, we need one more assumption on the
 $\lambda_i$, $i\in J_o$. We fix $i$ again. Let 
$$ {\mathscr I} = \{ e\in E_i^o\ ; \ \lambda_i^o (e{\mathfrak
o}_{E_i^o})\subset {\mathfrak p}_{F_o}\} \ .
$$
\noi This is an ${\mathfrak o}_{E_i^o}$-lattice in $E_i^o$ and must
have the form $t{\mathfrak p}_{E_i^o}$, for some $t\in
(E_i^o)^{\times}$. So replacing $\lambda_i$ by $e\mapsto \lambda_i
(tx)$, we may assume that ${\mathscr I}={\mathfrak p}_{E_i^o}$. In
the following  we assume that the linear forms $\lambda_i$,
$i\in J_o$, have this property.

\smallskip

\noi {\bf Lemma 5.4.} {\it Fix $i\in J_o$. Let $\lambda_i^{1}$,
$\lambda_i^{2}~:$ $E_i\rightarrow F$  be two linear forms as above and let
$h_i^{1}$, $h_i^{2}$ be the corresponding $\epsilon$-hermitian
forms on $V_i$ (i.e. $h_i^1$ and $h_i^{2}$ satisfy (5.3)). Then
there exists $u\in {\mathfrak o}_{E_i^{o}}^{\times}$ such that
$h_i^2 =uh_i^1$.}
\smallskip

\noi {\it Proof}. Since $h_i^1$ and $h_i^2$ induce the same
involution on ${\rm End}_{E_i}V_i$, there exists $u\in
E_i^{\times}$ such that $h_i^{2}=uh_i^{1}$. The fact that
$h_i^{1}$ and $h_i^{2}$ are both $\epsilon$-hermitian
with respect to $\sigma_{E_i}$ implies that $u$ lies in $E_i^{o}$.
Condition (5.3) writes
$$
h(v,w) =\lambda_i^{1}(h_i^{1}(v,w)) =
\lambda_i^{2}(uh_i^{1}(v,w))\ , \ v,w\in V_i \ .
$$
\noi So  $\lambda_i^{1}(e) = \lambda_i^{2}(ue)$, $e\in E_i$. By
applying ${\rm Tr}_{F/F_o}$ to this equality, we get
$\lambda_i^{o,1}(e) = \lambda_i^{o,2}(ue)$, $e\in
E_i^{o}$. Hence
$$
{\mathfrak p}_{E_i^{o}}=\{ e\in E_i^o\ ; \
\lambda_i^{o,1}(e{\mathfrak o}_{E_i^{o}})\subset {\mathfrak
p}_{F_o}\}
$$
$$
= \{ e\in E_i^o\ ; \ \lambda_i^{o,2}(ue{\mathfrak o}_{E_i^{o}}
\subset ) {\mathfrak p}_{F_o}\} = u^{-1}{\mathfrak p}_{E_i^{o}}\ .
$$
\noi So $u\in {\mathfrak o}_{E_i^{o}}^{\times}$ as required.

\smallskip

 Let us fix $i$. Let $L$ be an ${\mathfrak o}_{E_i^o}$-lattice in
 $V_i$. Then $L$ has a dual $L^{\sharp}$ relative to the form
 $h_{\vert V_i}$ and a dual $L^{\sharp_i}$ relative to the form $h_i$. 

\noi {\bf Lemma 5.5.} {\it The lattices $L^{\sharp}$ and
$L^{\sharp_i}$ coincide.} 
\smallskip

\noi {\it Proof}. We have 

\smallskip

\noi $L^{\sharp} =\{ v\in V_i \ ; \ h(v,L)\subset {\mathfrak p}_F\}$

 $= \{ v\in V_i\ ; \ {\rm Tr}_{F/F_o} h(v,L)\subset {\mathfrak
p}_{F_o}\}$

 $=\{ v\in V_i\ ; \ \lambda_{o}\circ {\rm Tr}_{E_i /E_i^o} h_i
(v,L)\subset {\mathfrak p}_{F_o}\}$

 $= \{ v\in V_i \ ; \  {\rm Tr}_{E_i /E_i^o} h_i
(v,L)\subset {\mathfrak p}_{E_i^{o}}\}$

 $= \{ v\in V_i \ ; \ f(v,L)\subset {\mathfrak p}_{E_i}\}$

 $=L^{\sharp_i}$,

\smallskip

\noi where the second and fifth equalities hold because $F/F_o $ and
$E_i /E_i^{o}$ are at worst tamely ramified.

\bigskip

\noi {\bf 6. Embedding the building of the centralizer}
\bigskip

 We keep the notation as in the previous section. Assume for a moment
 that the extensions $E_i /F$, $i\in J_o\cup J$, are separable.  Then the
 group $H$ is naturally the group of rational points of a reductive
 $F$-group $\boldsymbol H$. Indeed each $H_i$, $i\in J_o\cup J$, is
 naturally the group of rational points of a classical $E_i$-group
 ${\boldsymbol H}_i$ (we do not need $E_i /F$-separable here) and 
$$
{\boldsymbol H}\simeq \prod_{i\in J_o\cup J_+}{\rm Res}_{E_i /F} 
{\boldsymbol H}_i\ .
$$
\noi The (enlarged) affine building of $\boldsymbol H$, $I^1_{\beta }:= 
I^1({\boldsymbol H},F)$, is the cartesian product of the (enlarged) affine 
buildings $I^1({\rm Res}_{E_i /F} {\boldsymbol H}_i ,F)$, $i\in J_o\cup J_+$. 
For all $i$, the (enlarged) buildings $I^1({\rm Res}_{E_i /F} 
{\boldsymbol H}_i ,F)$ and $I^1({\boldsymbol H}_i ,E_i)$ identify 
canonically. Note also that, for $i\in J_o$, the centre of ${\boldsymbol H}_i$ 
is compact so the enlarged building is also the non-enlarged building; in 
particular, if $J=\emptyset$ then all the buildings involved are non-enlarged.

Since we do not want any restriction on the extensions $E_i /F$, we shall take 
as a definition of the (enlarged) building $I^1_{\beta}$ attached to
the group $H$:
$$
I^1_{\beta} := \prod_{i\in J_o\cup J_+} I^1({\boldsymbol H}_i ,E_i)
\leqno \hbox{ (6.1)}
$$
\noi We abbreviate $I^1_i =I^1({\boldsymbol H}_i ,E_i )$, 
$i\in J_o\cup J_+$.
\smallskip

 We are going to construct a map $j_{\beta}~:$ $I^1_{\beta} \rightarrow
I$.  
 We normalize the lattice-functions in ${\rm Latt}^{1}_{{\mathfrak
 o}_{E_i}}(V_i) $ by $\Lambda_i (r+ v_i (\pi_i ))={\mathfrak
 p}_{E_i}\Lambda_i (r)$, $r\in \mathbb R$, where, for each $i$, $\pi_i$
 denotes a uniformizer of $E_i$ and $v_i$ the unique extension of $v$
 to a valuation of $E_i$. It is straightforward that we have a well
 defined map
$$
{\tilde j}_{\beta}~:\ \prod_{i\in J_o\cup J} {\rm Latt}^{1}_{{\mathfrak
 o}_{E_i}}(V_i ) \longrightarrow {\rm Latt}^{1}(V)
$$
$$
\left(\Lambda_i\right)_{i\in J_o\cup J} \mapsto \bigoplus_{i\in J_o\cup J} 
\Lambda_i
$$
\noi where $\left(\bigoplus_{i\in J_o\cup J} \Lambda_i\right)(r) = 
\bigoplus_{i\in J_o\cup J} \Lambda_i(r)$, for $r\in \mathbb R$. This map is
clearly injective and equivariant for the action of group 
$\displaystyle \prod_{i\in J_o\cup J}{\rm Aut}_{E_i} V_i \subset {\rm
Aut}_F V$.  

For $i\in J_o$, we denote by $\sharp_i$ the involution on ${\rm
Latt}^{1}_{{\mathfrak o}_{E_i}} (V_i )$ attached to $h_i$, and by ${\rm
Latt}^{1}_{{\mathfrak o}_{E_i},h_i} (V_i ) \subset {\rm
Latt}^{1}_{{\mathfrak o}_{E_i}} (V_i )$ the set of fixed points. 
For $i\in J$, we denote be $\sharp_i$ the map ${\rm
Latt}^{1}_{{\mathfrak o}_{E_i}} (V_i )\to {\rm
Latt}^{1}_{{\mathfrak o}_{E_{-i}}} (V_{-i})$ given by
$$
\Lambda_i^{\sharp_i} (r) = \{v\in V_{-i}\ ; \ h(v,\Lambda_i(-r+))\subset 
\pf\} \ .
$$
for $\Lambda_i\in {\rm Latt}^{1}_{{\mathfrak o}_{E_i}} (V_i )$.

We define an involution $b$ on $\displaystyle \prod_{i\in J_o\cup J}{\rm
Latt}^{1}_{{\mathfrak o}_{E_i}}(V_i )$ by
$$
\left(\Lambda_i\right)_{i\in J_o\cup J}^b = 
\left(\Lambda_{-i}^{\sharp_{-i}}\right)_{i\in J_o\cup J}\ ,
$$
Then we have a bijection 
$$
\iota_h: \prod_{i\in J_o}{\rm Latt}^{1}_{{\mathfrak o}_{E_i},h_i}(V_i ) 
\times\prod_{i\in J_+}{\rm Latt}^{1}_{{\mathfrak o}_{E_i}}(V_i ) \to 
\left(\prod_{i\in J_o\cup J} {\rm Latt}^{1}_{{\mathfrak o}_{E_i}}(V_i )\right)^b,
$$
given by $\left(\Lambda_i\right)_{i\in J_o\cup J_+} \mapsto 
\left(\Lambda_i\right)_{i\in J_o\cup J}$, with 
$\Lambda_{-i}=\Lambda_i^{\sharp_i}$, for $i\in J_+$.

\smallskip

\noi {\bf Lemma 6.2.} {\it For $x\in \displaystyle\prod_{i\in J_o\cup J}{\rm
Latt}^{1}_{{\mathfrak o}_{E_i}}(V_i )$, we have ${\tilde
j}_{\beta}(x^b ) ={\tilde j}_{\beta}(x)^{\sharp_h}$. In particular
${\tilde j}_{\beta}\circ\iota_h$ maps  $\displaystyle\prod_{i\in J_o}{\rm
Latt}^{1}_{{\mathfrak o}_{E_i},h_i}(V_i )\times\prod_{i\in J_+}{\rm
Latt}^{1}_{{\mathfrak o}_{E_i}}(V_i)$ into  ${\rm Latt}^{1}_{h}(V)$.}
\smallskip

\noi {\it Proof}. Fix $\left(\Lambda_i\right)_{i\in J_o\cup J}\in 
\displaystyle\prod_{i\in J_o\cup J} {\rm Latt}^{1}_{{\mathfrak o}_{E_i}}V_i$ 
and set $\Lambda ={\tilde j}_{\beta}
\left(\left(\Lambda_i\right)_{i\in J_o\cup J_+}\right)$. We have
$$
\Lambda^{\sharp_h} (r)=\Lambda (-r+)^{\sharp_h} =\{ v\in V\ ; \ 
h(v,\Lambda (-r+))\subset {\mathfrak p}_F\}\ ,\ r\in {\mathbb R}\ .
$$
\noi Fix $r\in\mathbb R$. We have 
$$ 
\Lambda (-r+ ) =\bigoplus_{i\in J_o\cup J}\Lambda_i (-r+ ) \ .
$$
\noi Let $v=\sum_{i\in J_o\cup J}v_i$, with $v_i\in V_i$, be an element of
$V$. Since $V_i^\perp=\bigoplus_{k\ne -i} V_k$, we have 
$v\in \Lambda^{\sharp_h}(r)$ if and only if $h(v_{-i},\Lambda_{i} (-r+ )) 
\subset {\mathfrak p}_F$, for all $i$, that is if $v_{-i} \in 
\Lambda_{i}^{\sharp_i}(r)$, for all $i$ (by lemma (5.5) for $i\in J_o$ or by 
definition for $i\in J$); the lemma follows.
\smallskip

With the notation of {\S}4, for each set $\{j_i\}_{i\in J_+}$ of maps
$j_i: I_i^1\to{\rm Latt}^{1}_{{\mathfrak o}_{E_i}}(V_i)$ given by 
Proposition 4.3, we define a map $\displaystyle j_{\beta}~:\ 
\prod_{i\in J_o\cup J_+}I^1_i \rightarrow I$ by 
$$
j_{\beta} = j_{h}^{-1}\circ {\tilde j}_{\beta}\circ\iota_h\circ
\left(\prod_{i\in J_o} j_{h_i}\times \prod_{i\in J_+} j_{i}\right) \ .
$$
These maps depend {\it a priori} on the forms $h$, and $h_j$, $j\in J_o$.
\smallskip

\noi {\bf Theorem 6.3}. {\it Each map $j_{\beta}$ is injective,
$H$-equivariant. The set of such maps (as $\{j_i\}_{i\in J_+}$ varies) 
depends only on the involution $\sigma$.}
\medskip

\noi In particular, if $J=\emptyset$ then there is a unique map $j_\beta$, 
depending only on the involution $\sigma$. 

\medskip

\noi {\it Proofs}. The first two properties are
straightforward. Assume that $h' =uh$, $u\in F^{\times}$, is another
$\epsilon$-hermitian form on $V$, with respect to $\sigma_F$, defining the same
involution $\sigma$ on $\tGfr$. Then $u\in F_o$. For $i\in J_o$,
let $h_{i}'$ be an $\epsilon$-hermitian form on $V_i$ satisfying 
$$
uh(v,w)= \lambda_{i}' (h_{i}' (v,w))\, \ v,w\in V_i\ ,
$$
\noi where  the $\lambda_{i}'~:$ $E_i \rightarrow F$ are linear forms
as above. Then by lemma (5.4), for all $i\in J_o$, there exists $u'_i\in
{\mathfrak o}_{E_{i}^{o}}^{\times}$ such that $u^{-1}h_{i}' =u'_i h_i$, that is
$h_{i}' = uu'_i h_{i}$. 

\noi Let $\{j_i\}_{i\in J_+}$ be as above; we show that, for a suitable 
choice of $\{j'_i\}_{i\in J_+}$, we have
$$
j_{h}^{-1}\circ {\tilde j}_{\beta}\circ\iota_h\circ j_1
= j_{h'}\circ^{-1} {\tilde j}_{\beta}\circ\iota_{h'}\circ j_1',
$$
and the result follows. 

\smallskip
\noi By Lemma (3.6), for $i\in J_+$, for all $x_i \in I^1_i$, we have 
$j_{h_{i}'}(x_i )= j_{h_i}(x_i )-v(uu'_i)/2 = j_{h_i}(x_i )-v(u)/2$. 
For $i\in J_+$, we choose $j_i'$ such that $j_i'(x)=j_i(x)-v(u)/2$ for 
$x\in I_i^1$, that is $j_i'\circ j_i^{-1}(\Lambda_i)=\Lambda_i-v(u)/2$ for 
$\Lambda_i\in{\rm Latt}^{1}_{{\mathfrak o}_{E_i}}(V_i)$. We abbreviate
$$
j=\prod_{i\in J_o} j_{h_i}\times \prod_{i\in J_+} j_{i}, \qquad
j'=\prod_{i\in J_o} j_{h'_i}\times \prod_{i\in J_+} j'_{i};
$$
then, for $\left(\Lambda_i\right)_{i\in J_o\cup J_+}\in
\displaystyle\prod_{i\in J_o}{\rm Latt}^{1}_{{\mathfrak o}_{E_i},h_i}(V_i ) 
\times\prod_{i\in J_+}{\rm Latt}^{1}_{{\mathfrak o}_{E_i}}(V_i )$, we have
$$
j'\circ j^{-1} \left( \left(\Lambda_i\right)_{i\in J_o\cup J_+}\right) = 
\left(\Lambda_i -v(u)/2\right)_{i\in J_o\cup J_+}.
$$

\medskip
\noi It is also straightforward to check that
$$
\iota_{h'}\left(\left(\Lambda_i-v(u)/2\right)_{i\in J_o\cup J_+}\right)=
\iota_h\left(\left(\Lambda_i\right)_{i\in J_o\cup J_+}\right)-v(u)/2,
$$
for $\left(\Lambda_i\right)_{i\in J_o\cup J_+}\in
\displaystyle\prod_{i\in J_o}{\rm Latt}^{1}_{{\mathfrak o}_{E_i},h_i}(V_i ) 
\times\prod_{i\in J_+}{\rm Latt}^{1}_{{\mathfrak o}_{E_i}}(V_i )$. 
Then we have
\begin{eqnarray*}
{\tilde j}_{\beta}\circ\iota_{h'}\circ j'\circ j^{-1} 
\left( \left(\Lambda_i\right)_{i\in J_o\cup J_+}\right) 
&=& {\tilde j}_\beta\circ\iota_{h'}
\left(\left(\Lambda_i-v(u)/2\right)_{i\in J_o\cup J_+}\right) \\ 
&=& {\tilde j}_\beta\left(\iota_h\left(
\left(\Lambda_i\right)_{i\in J_o\cup J_+}\right)-v(u)/2\right) \\
&=&{\tilde j}_\beta\circ\iota_h
\left(\left(\Lambda_i\right)_{i\in J_o\cup J_+}\right)-v(u)/2.
\end{eqnarray*}

\noi By Lemma (3.6) again, we have $j_{h'}(x)=j_{h}(x)-v(u)/2$, $x\in
I$, that is $\Lambda -v(u)/2 = j_{h'}\circ j_{h}^{-1}(\Lambda)$,
$\Lambda\in {\rm Latt}^{1}_{h}(V)$. So
$$
j_{h'}\circ j_{h}^{-1}\circ {\tilde j}_{\beta}\circ \iota_h
 ={\tilde j}_{\beta}\circ\iota_{h'}\circ j'\circ j^{-1}\ ,
$$
\noi and the lemma follows.

\bigskip

\noi {\bf 7. Affine structures}
\medskip
  
We keep the notation as in the previous sections. For 
$x=(x_i)_{i\in J_o\cup J_+}$, $y=(y_i)_{i\in J_o\cup J_+}$ in $I^1_{\beta} =
\displaystyle \prod_{i\in J_o\cup J_+}I^1_i$ and $t\in [0,1]$, we define the 
barycenter $tx+(1-t)y$ to be 
$$
(tx_i +(1-t)y_i)_{i\in J_o\cup J_+} \ .
$$
\noi  We define the barycenter of two points in $\displaystyle 
\prod_{i\in J_o\cup J_+}{\rm Latt}^{1}_{{\mathfrak o}_{E_i}}(V_i )$ in a 
similar way. Since, for $i\in J_o$,  
${\rm Latt}^{1}_{{\mathfrak o}_{E_i} ,h_i}(V_i )$ is convex in 
${\rm Latt}^{1}_{{\mathfrak o}_{E_i}}(V_i )$, the subset
$\displaystyle\prod_{i\in J_o}{\rm Latt}^{1}_{{\mathfrak o}_{E_i},h_i}(V_i ) 
\times\prod_{i\in J_+}{\rm Latt}^{1}_{{\mathfrak o}_{E_i}}(V_i )$
of $\displaystyle \prod_{i\in J_o\cup J_+}{\rm Latt}^{1}_{{\mathfrak
o}_{E_i}}(V_i )$ is convex also. 

\medskip

\noi {\bf Proposition 7.1.} {\it Let $\beta$ be as in {\S}5. Then each
map $j_{\beta}$ is affine: for all $x$, $y\in I_{\beta}^{1}$, $t\in
[0,1]$, we have
$$
j_{\beta}(tx +(1-t)y) = tj_{\beta}(x)+(1-t)j_{\beta}(y)\ .
$$}

\noi {\it Proof}. By construction it suffices to prove that the maps  
${\tilde j}_{\beta}$
and $\iota_h$ are affine. We begin with ${\tilde j}_{\beta}$. 
Let $(\Lambda_i)_{i\in J_o\cup J}$, $(M_i)_{i\in J_o\cup J}$ 
be elements of $\displaystyle
\prod_{i\in J_o\cup J}{\rm Latt}^{1}_{{\mathfrak o}_{E_i}}(V_i)$. We must prove 
that 
$$
\bigoplus_{i\in J_o\cup J} (t\Lambda_i + (1-t)M_{i})
= t\left(\bigoplus_{i\in J_o\cup J}\Lambda_i\right)+
(1-t)\left(\bigoplus_{i\in J_o\cup J}M_i\right).
$$
\noi Let us recall the construction of the barycenter of two lattice
functions (we do it for ${\rm Latt}^{1}(V)$). Let $\Lambda$, $M\in
{\rm Latt}^{1}(V)$. There exists an $F$-basis $(e_1 ,\dots ,e_n )$
of $V$ which splits both $\Lambda$ and $M$~: there exist constants
$\lambda_1 ,\dots ,\lambda_n $, $\mu_1 ,\dots ,\mu_n $ in $\mathbb R$
such that 
$$
\Lambda (r) = \bigoplus_{k=1,\dots ,n}{\mathfrak p}_{F}^{\lceil
r+\lambda_k \rceil}e_k\ ,\ 
M (r) = \bigoplus_{k=1,\dots ,n}{\mathfrak p}_{F}^{\lceil
r+\mu_k \rceil}e_k\ , \ r\in {\mathbb R}\ .
$$
\noi  Then for $t\in [0,1]$, $t\Lambda +(1-t)M$ is given by 
$$
(t\Lambda +(1-t)M)(r) = \bigoplus_{k=1,\dots ,n}{\mathfrak p}_{F}^{\lceil
r+t\lambda_k +(1-t)\mu_k \rceil}e_k\ ,\ r\in {\mathbb R}\ .
$$

\noi The proof that ${\tilde j}_{\beta}$ is affine is then to construct a 
common splitting basis for $\bigoplus_{i\in J_o\cup J}\Lambda_i$ and
$\bigoplus_{i\in J_o\cup J}M_i$ from bases ${\mathscr B}_i$ of $V_i$,
$i\in J_o\cup J$, where ${\mathscr B}_i$ splits $\Lambda_i$ and
$M_i$. We leave this easy exercise to the reader.

\medskip

\noi Now we turn to $\iota_h$. Suppose $i\in J_+$ and 
$\Lambda_i\in{\rm Latt}^{1}_{{\mathfrak o}_{E_i}}(V_i )$, and let 
$(e_1 ,\dots ,e_n)$ be an $E_i$-basis of $V_i$ which splits $\Lambda_i$. Let 
$(e_{-1},\dots,e_{-n})$ be the dual $E_{-i}$-basis of $V_{-i}$, such that 
$h(e_{-k},e_l)=\delta_{kl}$, for $1\le k,l\le n$. It is straightforward to 
check that this basis splits $\Lambda_i^{\sharp_i}$ and that, 
$$
\hbox{if}\quad\Lambda_i (r) = \bigoplus_{k=1,\dots ,n}
{\mathfrak p}_{E_i}^{\lceil r+\lambda_k \rceil}e_k\quad\hbox{then}\quad  
\Lambda_i^{\sharp_i}(r) = \bigoplus_{k=1,\dots ,n}
{\mathfrak p}_{E_{-i}}^{\lceil r-\lambda_k \rceil}e_{-k}.
\leqno{(7.2)}
$$

\noi To show that $\iota_h$ is affine, we just need to check that, for 
$i\in J_+$, $\Lambda_i, M_i\in{\rm Latt}^{1}_{{\mathfrak o}_{E_i}}(V_i )$ and 
$t\in [0,1]$, we have
$$
\left(t\Lambda_i+(1-t)M_i\right)^{\sharp_i}=
t\Lambda_i^{\sharp_i}+(1-t)M_i^{\sharp_i}.
$$
The details of the proof -- which is to choose an $E_i$-basis of $V_i$ 
which splits both $\Lambda_i$ and $M_i$, take its dual basis and then 
use (7.2) -- are again left to the reader.
\bigskip

\noi {\bf 8. The image of an apartment}
\medskip

We keep the notation of the previous sections. We will show that the
image of an apartment of $I_\beta^1$ under each map $j_\beta$ is
contained in an apartment of $I$. 

Given a Witt decomposition $V=V_+\oplus V_o\oplus V_-$, with basis
$(e_l)_{l=1,...,r}$ of $V_+$ and the dual basis $(e_{-l})_{l=1,...,r}$
of $V_-$ (as in \S2), we get a (self-dual) decomposition 
$$
V=\bigoplus_{l=1}^{r} V^l\oplus V_o\oplus\bigoplus_{l=1}^{r} V^{-l},
$$
where $V^l=Fe_l=\left(\bigoplus_{k\ne -l}
V^l\oplus V_o\right)^\perp$. Such a decomposition (which we will also 
call a Witt decomposition) corresponds to the choice of an apartment 
$\mathscr A$ in $I$: in
terms of lattice functions, $j_h(\mathscr A)$ is the set of self-dual
lattice functions $\Lambda$ such that 
$$
\Lambda(s)=\bigoplus_{l=1}^{r} (V^l\cap\Lambda(s))\oplus 
(V_o\cap\Lambda(s))\oplus\bigoplus_{l=1}^{r}
(V^{-l}\cap\Lambda(s)),\qquad\hbox{for all }s\in\mathbb R,
$$
that is, $\Lambda$ is {\it split\/} by the decomposition (cf.\
Proposition 4.1).

Similarly, the choice of an (enlarged) apartment $\mathscr A^1$ in
$I^1_{\beta} = \displaystyle \prod_{i\in J_o\cup J_+}I^1_i$ is given
by similar $E_i$-decompositions of $V_i$ for $i\in J_o$ and (without the
self-duality restriction) $i\in J_+$.

\medskip

\noi {\bf Proposition 8.1.} {\it Let $\mathscr A^1$ be an (enlarged) 
apartment of $I^1_\beta$. Then there is an apartment $\mathscr A$ of $I$ 
such that $j_\beta(\mathscr A^1)\subset\mathscr A$.}

\medskip

\noi {\it Proof}. We write $\mathscr A^1=\prod_{i\in J_o\cup
J_+}\mathscr A^1_i$, with $\mathscr A^1_i$ an (enlarged) apartment in
$I^1_i$. 

As above, for each $i\in J_o$, the apartment $\mathscr A^1_i$
corresponds to a Witt $E_i$-decomposition of $V^i$
$$
V_i=\bigoplus_{l=1}^{r_i} V_i^l\oplus V_{i,o}\oplus\bigoplus_{l=1}^{r_i}
V_i^{-l},
$$ 
with $V_i^l=\left(\bigoplus_{k\ne -l} V_i^l\oplus V_{i,o}\right)^\perp$,
$\dim_{E_i}V_i^l=1$ and $r_i$ the ($E_i$-)Witt index of $V_i$. We
write ${\rm Latt}_{\mathfrak o_{E_i}}^{\mathscr A^1}(V_i)$ for the set of
lattice functions split by this decomposition, and ${\rm Latt}_{\mathfrak
o_{E_i}, h_i}^{\mathscr A^1}(V_i)$ for the subset of self-dual 
lattice functions, so that $j_{h_i}(\mathscr A^1_i)={\rm Latt}_{\mathfrak
o_{E_i}, h_i}^{\mathscr A^1}(V_i)$.

Also, for each $i\in J_+$, the apartment $\mathscr A^1_i$
corresponds to a decomposition of $V_i$ as a sum of $1$-dimensional
$E_i$-subspaces,
$$
V_i=\bigoplus_{l=1}^{r_i} V_i^l,
$$
with $r_i=\dim_{E_i}V_i$. As above, $j_{i}(\mathscr A^1_i)={\rm Latt}_{\mathfrak
o_{E_i}}^{\mathscr A^1}(V_i)$, the set of lattice functions split by
this decomposition.

We also take the dual splitting of $V_{-i}$ as a sum of $1$-dimensional
$E_{-i}$-subspaces, 
$$
V_{-i}^{l}=\left(\bigoplus_{k\ne l} V_i^k\right)^\perp.
$$
We remark that, if $\Lambda\in {\rm Latt}_{\mathfrak
o_{E_i}}^{\mathscr A^1}(V_i)$ then $\Lambda_i^{\#_i}$ is split by this
decomposition. 

\medskip

Now, for $i\in J_o\cup J_+$ and $1\le l\le r_i$, we decompose $V_i^l$
as a sum of $1$-dimensional $F$-subspaces as follows: fix $v\in
V_i^l$, $v\ne 0$, and let
${\mathscr B}_i$ be an $F$-basis for $E_i$ which splits the
$\mathfrak o_F$-lattice sequence $s\mapsto\mathfrak p_{E_i}^{\lceil
s/e(E_i/F)\rceil}$; then we take the decomposition
$$
V_i^l=\bigoplus_{b\in{\mathscr B}_i} Fbv.
$$
Note that any $\mathfrak o_{E_i}$-lattice sequence in $V_i^l$ is split
by this decomposition.
For $i\in J_o$, we also take the dual decomposition of $V_i^{-l}$ and,
for $i\in J_+$, the dual decomposition of $V_{-i}^l$. 

\medskip

Now we need 
to decompose the anisotropic parts $W:=\oplus_{i\in J_o} V_{i,o}$ suitably, 
for which we cheat.
Let $\boldsymbol G_o$ denote the classical group associated to the 
restriction of the form $h$ to $W$ and, for $i\in J_o$, let 
$\boldsymbol H_{i,o}$ denote the group associated to the restriction of 
the form $h_i$ to $V_{i,o}$. Note that the groups $H_{i,o}$ are compact 
so the building $I^1_{\beta,o}:=I^1({\boldsymbol H}_{i,o} ,E_i)$ is reduced 
to a point. 

Now, our constructions in \S6 give an embedding of $I^1_{\beta,o}$ in 
the building $I^1_o:=I^1(\boldsymbol G_o,F)$ and the image is 
certainly contained in some apartment. Hence there is a Witt 
$F$-decomposition of $W$ which splits the (unique) self-dual 
lattice sequence in $W$ corresponding to $I^1_{\beta,o}$, and this is 
the decomposition we take.

\medskip

Altogether, we have described a Witt $F$-decomposition of $V$, 
which corresponds to an apartment $\mathscr A$ of $I$. We denote 
by ${\rm Latt}_{\mathfrak o_F, h}^{\mathscr A}(V)$ the set of 
self-dual lattice functions in $V$ which are split by this splitting, 
so that $j_{h}(\mathscr A)={\rm Latt}_{\mathfrak o_F, h}^{\mathscr A}(V)$.

Finally, by construction it is clear that $\tilde j_\beta\circ\iota_h$ maps 
$\displaystyle
\prod_{i\in J_o} {\rm Latt}_{\mathfrak o_{E_i}, h_i}^{\mathscr A^1}(V_i)
\times\prod_{i\in J_+} {\rm Latt}_{\mathfrak o_{E_i}}^{\mathscr A^1}(V_i)$ 
into ${\rm Latt}_{\mathfrak o_F, h}^{\mathscr A}(V)$ so 
$j_\beta(\mathscr A^1)\subset\mathscr A$, as required.

\bigskip

\noi {\bf 9. Compatibility with Lie algebra filtrations}
\medskip

 In this section, we fix $H_k$-equivariant identifications $j_k$~:
 $I^{1}(H_k ,E_k)\rightarrow {\rm Latt}^{1}_{\ofr_{E_k}}( V_k )$, $k\in
 J^{+}$. They give rise to the map $j_{\beta}$~:
 $I_{\beta}^{1}\rightarrow I(G,H)$ defined in {\S}6.
\smallskip

 Let $x\in I(G,F)=I^{1}(G,F)$, that we see as a self-dual lattice
 function $\Lambda$ in ${\rm Latt}_{h}^{1}(V)$. To $x$ we can
 associate a filtration $(\Gfr_{x,r})_{r\in \RR}$ of the Lie algebra
 $\Gfr$ as follows. First $x$ defines a filtration
 $(\tGfr_{x,r})_{r\in \RR}$ of $\tGfr$ by
$$
\tGfr_{x,r}=\{a\in \tGfr\ ; \ a\Lambda (s)\subset \Lambda (s+r),\ s\in
\RR\} , \ r\in \RR \ .
$$
\noi We then define
$$
\Gfr_{x,r}:=\tGfr_{x,r}\cap \Gfr =\{a\in \Gfr\ ; \ a\Lambda (s)\subset \Lambda (s+r),\ s\in
\RR\} , \ r\in \RR \ . \eqno \hbox{(1)}
$$
\noi Similarly a point $x$ of  $I_{\beta}^{1}$ defines a filtration
$(\Hfr_{x,r})_{r\in \RR}$ of $\Hfr$. Write $x=(x_k )_{k\in J\cup
J_o}$, $x_k \in I^{1}(H_k ,E_k )$; each $x_k$ corresponding to a
lattice function $\Lambda_k$ of ${\rm Latt}_{\ofr_{E_k}}(V_k )$ (with
$\Lambda_{k}^{\sharp_{k}}=\Lambda_{-k}$, $k\in J\cup J_o$). We then
define
$$
\Hfr_{x,r}:=\bigoplus_{k\in J^{+}\cup J_o}\Hfr_{x_k ,r}^{k} , \ r\in
\RR , \leqno \hbox{(2)}
$$
\noi where 
$$
\Hfr_{x_k ,r}^{k}=\{ a\in {\rm Lie}(H_k )\ ; \ a\Lambda_{k}(s) \subset
\Lambda_{k}(s+r), \ s\in \RR\},\ r\in \RR , \ k\in J^{+}\cup J_o\ .
$$

 The filtration $(\Hfr_{x,r})_{r\in \RR}$ only depends on the image
 $\bar x$ of $x$ in the non-enlarged building $I_{\beta}$. One can
 prove that for $x\in I(G,F)$, $(\Gfr_{x,r})_{r\in \RR}$ is the
 filtration of $\Gfr$ attached to $x$ defined by Moy and Prasad
 \cite{MP}. Similarly, when $\beta$ is semisimple and $x\in
 I^{1}(H,F)$,  $(\Hfr_{x,r})_{r\in \RR}$ is the filtration of $\Hfr$
 attached to $\bar x$ defined in loc. cit. The proof of this fact is
 announced by B. Lemaire and J.-K. Yu [BY].
\medskip

\noi {\bf Lemma 9.1}. {\it Let us see $\Hfr$ as being canonically
embedded in $\displaystyle \tHfr ={\rm End}_E V=\bigoplus_{k\in J\cup
J_o}{\rm End}_{E_k}V_k$ via
$$
(a_k)_{k\in J^{+}\cup J_o}\mapsto (b_k)_{k\in J\cap J_o}\ ,
$$
\noi where $b_k = a_k$, $k\in J_o$, and $b_{-k}=-a_{k}^{\sigma}$, $k\in
J^{+}$. Fix $x\in I^{1}_{\beta}$ as before and consider the
$\of$-lattice function in $V$ given by 
$$
\Lambda =\bigoplus_{k\in J\cup J_o}\Lambda_k \ {\rm (notation\  of \
{\S}6)}\ .
$$
For $r\in \RR$, let 
$$
\tHfr_{x,r}=\{ a\in \tHfr \ ; \ a\Lambda (s)\subset \Lambda (s+r),\
s\in \RR\} , \ r\in \RR\ .
$$
\noi Then we have $\Hfr_{x,r}=\tHfr_{x,r}\cap \Hfr$,  $r\in \RR$.}
\smallskip

\noi {\it Proof}. Indeed, for all $a=(a_k )_{k\in J\cup J_o}\in {\rm
End}_E V$, we have $a\in \tHfr_{x,r}\cap \Hfr$ if and only if
$a+a^{\sigma}=0$ and $a\Lambda (s)\subset \Lambda (s+r)$, $s\in \RR$,
i.e. 
$$
a_k \Lambda_{k}(s)\subset \Lambda_{k}(s+r) , \ s\in \RR , k\in J\cup
J_o\ .
$$
\noi For $k\in J_o$, these conditions can be rewritten $a_k\in {\rm
Lie}(H_k )$ and $a_k \Lambda_k (s)\subset \Lambda_k (s+r)$, $s\in
\RR$, that is $a_k \in \Hfr_{x,r}^{k}$, as required. For $k\in J$,
these conditions can be rewritten $a_{-k}=-a_{k}^{\sigma}$ and
$$
a_k \Lambda_{k}(s)\subset\Lambda_k (s+r) , \ s\in \RR \eqno \hbox{(a)}
$$
$$
-a_{k}^{\sigma}\Lambda_{k}^{\sharp_{k}}(s)\subset
 \Lambda_{k}^{\sharp_k}(s+r),\ s\in \RR \ .\eqno \hbox{(b)}
$$
\noi So we must prove that conditions (a) and (b) are equivalent. By
symmetry we only prove one implication. Applying the duality $\sharp_k$
on lattices of $V_k$ to inclusion (b), we obtain
$$
\Lambda_{k}((-s-r)+)\subset
[a_{k}^{\sigma}\Lambda_{k}^{\sharp_k}(s)]^{\sharp_k} ,\ s\in \RR , 
$$
\noi with
$$
[a_{k}^{\sigma}\Lambda_{k}^{\sharp_k}(s)]^{\sharp_k} =\{ v\in V_k \ ;
\ a_k v\in \Lambda_{k}((-s)+)\} , \ s\in \RR \ .
$$
\noi So we have
$$
a_{k}\Lambda_{k}((-s-r)+)\subset \Lambda_{k}((-s)+)\subset \Lambda_k
(-s) , \ s\in \RR\ ,
$$
\noi that is 
$$
a_{k}\Lambda (s+)\subset \Lambda_{k}(s+r), \ s\in \RR\ .
$$
\noi On each open interval $(u,v)$ where $\Lambda_k$ is constant, we
have 
$$
a_{k}\Lambda_k (s+)=a_k \Lambda_k (s) \subset \Lambda_{k}(s+r)\ ,
$$
and (a) is true for $s\in (u,v)$. Finally if $s_o$ is a jump of
$\Lambda_k$ with $\Lambda_k$ constant on $(t,s_o ]$, we have
$$
a_k \Lambda_k (s_o )=a_k \Lambda_k (s+)\subset \Lambda_k (s+r) , \
s\in (t,s_o)\ .
$$
\noi So
$$
a_k \Lambda_k (s_o )\subset \bigcap_{s\in (t,s_o )}\Lambda_{k}(s+r)
=\Lambda_k (s_o +r)\ ,
$$
\noi $\Lambda_k$ being left continuous, and (a) is then true for all
$s\in \RR$. 
\medskip

\noi {\bf Proposition 9.2.} {\it Let $x\in I_{\beta}^{1}$. Then we
have
$$ \Gfr_{j_{\beta}(x),r}\cap \Hfr =\Hfr_{x,r}, \ r\in \RR\ .
$$}
\smallskip

\noi {\it Proof}. Indeed, with the notation of (9.1) and by definition
of $j_{\beta}$, we easily see that   
$$
\tGfr_{j_{\beta}(x),r}\cap \tHfr = \tHfr_{x,r}\ .
$$
\noi So our result is now a corollary of (9.1) since $\Hfr =\Gfr \cap
\tHfr$. 
\bigskip

\noi {\bf 10.  A unicity result for the general linear group}
\medskip

 As in \cite{BL}{\S}I.2, we define an equivalence relation $\sim$ on ${\rm
 Latt}^{1}(V)$ by $\Lambda_1 \sim \Lambda_2$ if there exists $s\in
\RR$ such that $\Lambda_1 (s) =\Lambda_{2}(r+s)$, $s\in \RR$. Then
$\sim$ is compatible with the $\tilde G$-action and the quotient
${\rm Latt}_{\of}(V):= {\rm Latt}^{1}(V)/\sim$ is naturally a 
$\tilde G$-set. We shall denote by $\bar \Lambda$ an element of  ${\rm
 Latt}_{\of}(V)$, where $\Lambda$ is a representative in ${\rm
 Latt}^{1}(V)$. As a consequence of [BL]{\S}I.2 and [BT1], there
 is a unique affine and $\tilde G$-equivariant map $j~:$ ${\tilde
 I}\rightarrow {\rm Latt}_{\of}(V)$, where $\tilde I$ denotes the
 non-enlarged building of $\tilde G$. 

 We fix an element $\beta$ of $\tGfr$ satisfying

$$
E: F[\beta ]\ {\rm is \ a \ field}\ .
\leqno (\hbox{H})
$$
\noi As in {\S}5 we denote by $\tHfr ={\rm End_{E}}V$ the centralizer of
$\beta$ in $\tGfr$ and by ${\tilde H}={\rm Aut}_E V$ its centralizer
in $\tilde G$. There is a canonical identification of the non-enlarged
affine building ${\tilde I}_{\beta}$ of $\tilde H$ with the $\tilde
H$-set ${\rm Latt}_{\ofr_E}(V)$. Here we normalize the lattice functions
of ${\rm Latt}_{\ofr_{E}}^{1}(V)$ by the condition $\Lambda (s+v(\pi_{E}
))=\pi_E \Lambda (s)$, $s\in \RR$, where $\pi_E$ is a uniformizer of
$E$.

Any ${\bar \Lambda}\in {\rm Latt}_{\of}(V)$ defines a filtration
$(\tGfr_{{\bar \Lambda}, r})_{r\in \RR}$ by
$$
\tGfr_{{\bar \Lambda}, r}=\{ a\in {\rm End}_F V\ ; \ a\Lambda
(s)\subset \Lambda (r+s), \ s\in \RR \}\ .
$$
\noi Then the map ${\rm End}({\bar \Lambda})~:$ $r\mapsto\tGfr_{{\bar
\Lambda}, r}$ is an element of ${\rm Latt}^{1}\tGfr$. The map ${\bar
\Lambda}\mapsto {\rm End}({\bar \Lambda})$, ${\rm
Latt}_{\of}V\rightarrow {\rm Latt}^{1}\tGfr$  is a $\tilde
G$-equivariant injection (cf. \cite{BL}{\S}4) for the action of $G$ on
${\rm Latt}^{1}\tGfr$ by conjugation. Its image is ${\rm
Latt}^{2}\tGfr$.  From now on we shall canonically identify $\tilde I$
(resp. ${\tilde I}_{\beta}$  with
${\rm Latt}^2 \tHfr$). 

Let us recall the main result of \cite{BL}.
\smallskip

\noi {\bf Theorem 10.1}. {\it There exists a unique affine and $\tilde
H$-equivariant map ${\tilde j}_{\beta}~:$ ${\tilde
I}_{\beta}\rightarrow {\tilde I}$. It is injective, maps any apartment
into an apartment and is compatible with the Lie algebra filtrations in
the following sense:
$$
\tGfr_{{\tilde j}_{\beta},r}\cap \tHfr =\tHfr_{x,r},\ x\in {\tilde
I}_{\beta}, \ r\in \RR\ .
\leqno (\hbox{10.2})   
$$}
\smallskip

Let us recall how ${\tilde j}_{\beta}$ is constructed. If $x\in
{\tilde I}_{\beta}$ corresponds to ${\rm End}({\bar \Lambda})\in {\rm
Latt}^2 \tHfr$, then ${\tilde j}(x)$ simply corresponds to ${\rm
End}({\bar \Lambda})$, where $\Lambda$, an $\ofr_E$-lattice function
in $V$, is now considered as an $\of$-lattice function. 
\smallskip

\noi {\bf Theorem 10.3.} {\it Let $x\in {\tilde I}_{\beta}$ and $y\in
{\tilde I}$ satisfying
$$
\tGfr_{y,r}\cap \tHfr \supset \tHfr_{x,r},\ r\in \RR \ .
$$
\noi Then $y={\tilde j}_{\beta}(x)$. As a consequence the map ${\tilde
j}_{\beta}$ is characterized by property (10.2).}
\smallskip

\noi {\it Proof}.  Assume that $x$ and $y$ correspond to elements ${\bar
\Lambda}_x$ and ${\bar \Lambda}_{y}$ of ${\rm Latt}_{\ofr_E}(V)$ and
${\rm Latt}_{\ofr_F}(V)$ respectively. 
\smallskip

\noi {\bf Lemma 10.4.} {\it Under the assumption of (10.2),
$\Lambda_y$ is an $\ofr_E$-lattice function.}
\smallskip

\noi {\it Proof}. To prove that $\Lambda_y$ is an $\ofr_E$-lattice function
we must prove that it is normalized by $E^{\times}=\langle
\pi_E \rangle \ofr_{E}^{\times}$, or equivalently:
$$
x\tGfr_{y,r}x^{-1}= \tGfr_{y,r}, \ x\in E^{\times},\ r\in \RR\, .
\leqno (\hbox{10.5})
$$
\noi We first notice than $\ofr_{E}\subset
\tHfr_{x,0}\subset\tGfr_{y,0}$, so that $\ofr_{E}^{\times}\subset
\tGfr_{y,0}^{\times}$ and (10.5) is true for $x\in
\ofr_{E}^{\times}$. We are reduced to proving (10.5) when $x=\pi_E$. 

 We have $\pi_{E}\in \tHfr_{x,1/e}\subset \tGfr_{y,1/e}$ and
$\pi_{E}^{-1}\subset \tHfr_{x,-1/e}\subset\tGfr_{y,-1/e}$, where
$e=e(E/F)$. It follows that
$$
\pi_{E}\tGfr_{y,r}\pi_{E}^{-1}\subset
\tGfr_{y,1/e}\tGfr_{y,r}\tGfr_{y,-1/e}\subset \tGfr_{y,r},\ r\in \RR
\, .
\leqno (\hbox{10.6})
$$
\noi Consider the  duality ``$*$'' on subsets of $\tGfr$ given by 
$$
S^{*}=\{a\in \tGfr \ ; \ {\rm Tr}(aS)\subset \pfr_{F}\},\ S\subset
\tGfr ,
$$
\noi where $\rm Tr$ is the trace map. Recall ([BL](6.3)) that
$(\tGfr_{y,r})^{*}=\tGfr_{y, (-r)+}$, $r\in \RR$. Using a well known
property of the trace map, we observe that 
$$
(\pi_E \tGfr_{y,r}\pi_{E}^{-1})^{*}= \pi_E
(\tGfr_{y,r})^{*}\pi_{E}^{-1},\ r\in\RR\, .
$$
\noi So applying the duality to (10.6), we obtain
$$
\tGfr_{y,(-r)+}\subset \pi_{E}\tGfr_{y,(-r)+}\pi_{E}^{-1} , \ r\in
\RR\, .
$$
\noi We have proved that on each open interval $(r_1 ,r_2 )$ where the
lattice function $(\tGfr_{y,r})_{r\in \RR}$ is constant, we have both
containments
$$
\pi_{E}\tGfr_{y,r}\pi_{E}^{-1}\subset\tGfr_{y,r}\ {\rm and}\
 \pi_{E}\tGfr_{y,r}\pi_{E}^{-1}\subset\tGfr_{y,r}, \ r\in
\RR\, .$$
\noi So by continuity we have
$\pi_{E}\tGfr_{y,r}\pi_{E}^{-1}=\tGfr_{y,r}$, for all $r$, as
required.
\smallskip

Let us return to the proof of (10.3). Since $\Lambda_y$ is an
$\ofr_E$-lattice function, we have
$$
\tGfr_{y,r}\cap \tHfr =\tHfr_{x' ,r},\ r\in \RR ,
$$
\noi where $x'\in {\tilde I}_{\beta}$ is attached to ${\bar
\Lambda}_{y}$, $\Lambda_{y}$ being seen as an $\ofr_E$-lattice
function. So by injectivity of the map ${\rm
Latt}_{\ofr_E}^{1}(V)\rightarrow {\rm Latt}^2 \tHfr$, we have ${\bar
\Lambda}_x ={\bar \Lambda}_{y}$ and $y={\tilde j}_{\beta}(x)$ by definition.

\bigskip

\noi {\bf 11. A unicity result in the $1$-block case and a conjecture}  
\medskip

With the notation of {\S}5, we consider an element $\beta\in \Gfr$
satisfying:
$$
E:= F[\beta ]\subset \tGfr\ {\rm is\ a\ field\ and}\ \beta \not= 0\, .
\leqno (\hbox{11.1})
$$

We fix an $\epsilon$-hermitian form $h_E$ on the $E$-vector space $V$
relative to $\sigma_E$ and we assume that it satisfies (5.3) as well
as the condition ${\mathcal J}=\pfr_{E^o}$ of {\S}5. This allows us to
identify $I_{\beta}^{1}$ with ${\rm Latt}^{1}_{h_E}(V)$. Identifying
$I$ with ${\rm Latt}_{h}(V)$, the map $j_{\beta}$ of {\S}6 is simply
given by
$$
j_{\beta}(\Lambda )=\Lambda ,\ \Lambda\in {\rm
Latt}^{1}_{h_E}(V),$$
\noi where on the right hand side $\Lambda$ is considered
as an $\ofr_{F}$-lattice function.
\smallskip

\noi {\bf Theorem 11.2}. {\it Under the assumption (11.1), let $x\in
I_{\beta}^{1}$ and $y\in I$ satisfying
$$
\Gfr_{y,r}\cap \Hfr =\Hfr_{x,r},\ r\in \RR\, .
\leqno (\hbox{11.3})
$$
\noi Then $y=j_{\beta}(x)$. In particular the map $j_{\beta}$ is
characterized by  compatibility with the Lie algebra filtrations.}
\smallskip

\noi {\it Proof}. The point $x$ (resp. $y$) corresponds to a self-dual
lattice function $\Lambda_x \in {\rm Latt}^{1}_{h_{E}}(V)$
(resp. $\Lambda_y \in {\rm Latt}^{1}_{h}(V)$).  We may see $x$ and $y$
as points of ${\rm Latt}^{1}_{\ofr_E}(V)$ and ${\rm
Latt}^{1}_{\ofr_F}(V)$ respectively  and they give rise to filtrations of $\tHfr$
and $\tGfr$ as in {\S}9: $(\tHfr_{x,r})_{r\in \RR}$ and
$(\tGfr_{y,r})_{r\in \RR}$. Write
$$
\Gfr_{y,r}^{+}=\{a\in \tGfr_{y,r}\ ; \ a=a^{\sigma}\},\ r\in \RR
$$
\noi and
$$
\Hfr_{x,r}^{+}=\{a\in \tHfr_{x,r}\ ; \ a=a^{\sigma}\},\ r\in \RR
$$
\noi Since $2$ is invertible in $\ofr_F$, we have:
$$
\tGfr_{y,r}=\Gfr_{y,r}\oplus \Gfr_{y,r}^{+}\ {\rm and}\ 
\tHfr_{y,r}=\Hfr_{x,r}\oplus \Hfr_{x,r}^{+},\ r\in \RR\, .
$$
\noi Write
$$
r_o =v_{\Lambda_x}(\beta ):={\rm Sup}\{r\in \RR\ ; \ \beta\in
\tHfr_{x,r}\}\, .
$$
\noi Since $\beta\in E^{\times}$, it normalizes $\Lambda_x$ so that
$\beta\tHfr_{x,r}=\tHfr_{x,r+r_o}$, $r\in\RR$. Moreover since $\beta$
is central in $\tHfr$, we easily have that
$\Hfr_{x,r}^{+}=\beta\Hfr_{x,r-r_o}$, $r\in \RR$. Hence, for $r\in
\RR$, we have
$$
\Hfr_{x,r}^{+}=\beta (\Gfr_{y,r-r_o}\cap \Hfr )=\beta
(\Gfr_{y,r-r_o}\cap \tHfr)\subset \Gfr_{y,r}\cap \tHfr\, .
$$
\noi It follows that, for $x\in \RR$, we have:
$$
\tHfr_{x,r}=\Hfr_{x,r}\oplus\Hfr_{x,r}^{+}\subset
\Gfr_{y,r}\cap\tHfr\oplus \Gfr_{y,r}^{+}\cap \tHfr\subset
\tGfr_{y,r}\cap \tHfr\, .
$$
\noi By applying (10.3), we obtain ${\bar \Lambda}_{y}={\tilde
j}_{\beta}({\bar \Lambda}_x)$, that is ${\bar \Lambda}_{y}={\bar
\Lambda}_x$. In particular we have ${\rm End}(\Lambda_x )={\rm
End}(\Lambda_y )\in {\rm Latt}^{2}_{\sigma}\tHfr$. But by (3.5) we
have $\Lambda_x =\Lambda_y$, as required.
\medskip

 Let us give an example. Assume that $G={\rm Sp}_2 (F )={\rm SL}(2,F)$
 (here $F=F_o$) and take
 $\beta\in \Gfr$ such that $E/F$ is quadratic and ramified. Then $H$
 is the group $E^1$ of norm $1$ elements in $E$. The building of $H$
 is reduced to a point $\{ x\}$. The group $E^{\times}$ fixes a unique
 chamber $C$ of $I$ and $H\subset E^{\times}$ fixes $C$
 pointwise. There are infinitely many maps $j$ :
 $I_{\beta}^{1}\rightarrow I$ which are affine and $G$-equivariant;
 indeed $j(x)$ can be any point of $C$. On the other hand there is a
 unique map $j$ : $I_{\beta}^{1} \rightarrow I$ which is compatible
 with the Lie algebra filtrations: it maps $x$ to the isobarycenter of
 $C$.  
\medskip

 We conjecture that when $J=\emptyset$ (notation of {\S}5) then the
 map $j_{\beta}$  of {\S}6 is characterized by condition (11.3). We
 may  address the more general (but more informal)
 question. Being given two $F$-reductive groups $\boldsymbol H$ and
 $\boldsymbol G$, as well as a morphism of algebraic groups $\varphi$
 : ${\boldsymbol H}\rightarrow {\boldsymbol G}$, is there an
 affine and ${\boldsymbol H}(F)$-equivariant map $I({\boldsymbol
 H},F)\rightarrow I({\boldsymbol G},F)$ which is compatible with the Lie
 algebra filtrations defined by Moy and Prasad. When is it characterized by
 this last property?

\bigskip

\noi {\bf References}
  
\begin{itemize}

\bibitem[BK]{BK} C.J. Bushnell and P.C. Kutzko, {\it The admissible
dual of ${\rm GL}(N)$ via compact open subgroups}, Ann. of
Math. Studies 129, Princeton Univ. Press. 

\bibitem[BL]{BL} P. Broussous and B. Lemaire, {\it Building of ${\rm
GL}(m,D)$ and centralizers}, Transformations Groups, Vol.\ 7, No.\ 1,
2002, pp.~15--50.

\bibitem[BT1]{BT1} F. Bruhat and J. Tits, {\it Sch\'emas en groupes et 
immeubles des groupes classiques sur un corps local, 1\`ere partie: le groupe 
lin\'eaire g\'en\'eral}, Bull.\ Soc.\ Math.\ Fr., 112, 1984, pp.~259--301.

\bibitem[BT2]{BT2} F. Bruhat and J. Tits, {\it Sch\'emas en groupes et 
immeubles des groupes classiques sur un corps local, 2\`eme partie: groupes 
unitaires}, Bull.\ Soc.\ Math.\ Fr., 115, 1987, pp.~141--195.

\bibitem[BY]{BY} B. Lemaire and J.-K. Yu, Private communication.

\bibitem[Inv]{Inv} M.-A. Knus, A. Merkurjev, M. Rost and J.-P. Tignol, 
{\it The book of involutions}, AMS Colloquium publications, Vol.\ 44, 1998.

\bibitem[Mor]{Mor} L. Morris, {\it Some tamely ramified supercuspidal 
representations of symplectic groups}, Proc.\ London Math.\ Soc.\ (3), 
Vol.\ 63, 1991, pp.~519--551.

\bibitem[MP]{MP} A. Moy and G. Prasad, {\it Unrefined minimal
$K$-types for $p$-adic groups}, Invent. Math., Vol.\ 116, 1994, 
pp.~393--408.

\bibitem[Sch]{Sch} W. Scharlau, {\it Quadratic and Hermitian forms}, 
Grundlehren der Mathematischen Wissenschaften 270, Springer-Verlag, Berlin, 
1985.

\bibitem [S1]{S1} S. Stevens, {\it Semisimple strata for {$p$}-adic 
classical groups}, Ann.\ Sci.\ \'Ecole Norm.\ Sup.\ (4), Vol.\ 35, 
No.\ 3, 2002, pp.~423--435.

\bibitem [S2]{S2} S. Stevens, {\it Semisimple characters for {$p$}-adic 
classical groups}, Preprint, September 2003.
\end{itemize}
\bigskip

Universit\'e de Poitiers \hfill School of Mathematics

D\'epartement de Math\'ematiques \hfill University of East Anglia

UMR 6086 du CNRS \hfill Norwich NR4 7TJ

Bd Marie et Pierre Curie \hfill United Kingdom

T\'el\'eport 2 -- BP 30179 

86962 Futuroscope Chasseneuil Cedex 

France

\end{document}